\documentclass[14pt]{article}
\usepackage{amssymb}
\usepackage{amsmath}
\usepackage[dvips]{graphics}
\usepackage{epsfig}
\usepackage{xcolor}

\newtheorem{theorem}{Theorem}

\newtheorem{remark}[theorem]{Remark}

\begin{document}
\date{}
\author{ Aristides V. Doumas \\
National Technical University of Athens, Greece\\
and\\
Archimedes/Athena Research Center, Greece\\
adou@math.ntua.gr, aris.doumas@hotmail.com}

\title{An inequality involving alternating binomial sums}
\maketitle

\begin{abstract}
In this letter, we prove an inequality involving alternating binomial logarithmic sums by exploiting the variance of the logarithm of the maximum of independent and identically distributed exponential random variables. This inequality was introduced in our recent work [On the minimum of independent collecting processes via the Stirling numbers of the second kind, \textit{Statist. Probab. Lett.}, \textbf{185} (2022)].
\end{abstract}

\textbf{MSC 2020 Mathematics Classification.} 60G70; 33B15.

\section{Introduction}
The \textit{coupon collector problem (CCP)} concerns a population whose members belong to \(N\) different types (\textit{species}). For $j=1,2,\cdots,N$, let $p_{j}$ denote the probability that a member of the population is of type $j$, where $p_{j}>0$ and $\sum_{j=1}^{N}p_{j}=1$. We refer to the $p_{j}$'s as the coupon probabilities. The members of the population are sampled independently \textit{with replacement}, and their types are recorded. The classical form of the problem corresponds to the case where the \(N\) distinct items are \textit{uniformly} distributed, that is, when $p_{j}=1/N$.

Let us consider a multi-player version of the \textit{CCP}, in which $n$ players aim to complete a full set of $N$ different coupons that are \textit{uniformly} distributed. Let $T_{N(1)}, T_{N(2)},\cdots, T_{N(n)}$ be the random variables denoting the numbers of trials needed until all \(N\) types are detected (at least once) by each of the $n$ independent collectors. This variant of the \textit{CCP} was introduced in \cite{D}, where the main task was to study the variance of the random variable
\[
M_{N(n)} := \bigwedge_{i=1}^n T_{N(i)}, \quad  N \to \infty,
\]
that is, the minimum of the random variables $\left\{T_{N(1)}, T_{N(2)},\cdots, T_{N(n)}\right\}$. As shown in \cite{D},
\begin{equation}
V\left[M_{N(n)}\right]\sim \left(\frac{\pi^{2}}{6} +n w_{n}- n^{2}c^{2}_{n}\right)N^{2},\,\,\, N \rightarrow \infty, \label{pla}
\end{equation}
where
\begin{equation}
c_{n}:=-\frac{1}{n}\sum_{j=1}^{n}\left(-1\right)^{j}\binom{n}{j}\ln j \label{cp}
\end{equation}
and
\begin{equation}
w_{n}:=-\frac{1}{n}\sum_{j=1}^{n}\left(-1\right)^{j}\binom{n}{j}\left(\ln j\right)^{2}. \label{wp}
\end{equation}
As usual, the symbol $\sim$ means that if $a_{N}\sim b_{N}$, then $\lim_{N\rightarrow \infty}a_{N}/b_{N}\rightarrow 1$.

In this letter, we answer an open question posed in \cite{D}. More precisely, we prove that the coefficient of the leading asymptotic term of $V\left[M_{N(n)}\right]$ is strictly positive for every \textit{fixed} positive integer $n$, i.e.,
\begin{equation}
\frac{\pi^{2}}{6} > n^{2}c^{2}_{n}-n w_{n},\,\,\,n\in \mathbb{N}. \label{v}
\end{equation}
Inequality (\ref{v}) was posed as an open question in \cite{D}.

\section{A probabilistic approach}
We proceed as follows. Let $X_{1}, X_{2},\cdots, X_{n}$ be independent and identically distributed exponential random variables such that
\begin{equation*}
E[X_{1}]=1.
\end{equation*}
Set
\begin{equation*}
W:=\max \left(X_{1}, X_{2},\cdots, X_{n}\right).
\end{equation*}
It is well known, and easy to verify, that
\begin{equation*}
F_{W}(w)=\left(1-e^{-w}\right)^{n}\mathbf{1}_{(0, \infty)}(w).
\end{equation*}
It follows that
\begin{equation*}
f_{W}(w)=ne^{-w}\left(1-e^{-w}\right)^{n-1}\mathbf{1}_{(0, \infty)}(w),
\end{equation*}
where, as usual, $\mathbf{1}_A$ denotes the indicator function of the event $A$. Next, define
\begin{equation*}
Y:=\ln W.
\end{equation*}
Hence,
\begin{equation*}
f_{Y}(y)=ne^{y}e^{-e^{y}}\left(1-e^{-e^{y}}\right)^{n-1},\,y\in \mathbb{R}.
\end{equation*}

Our next task is to derive the first and second moments of the random variable $Y$. We have
\begin{align}
E[Y]=&n\int_{-\infty}^{\infty}y\,e^{y}e^{-e^{y}}\left(1-e^{-e^{y}}\right)^{n-1}dy \nonumber\\
=&n\int_{0}^{\infty}\ln x\,\, e^{-x}\left(1-e^{-x}\right)^{n-1}dx  \nonumber\\
=&n\sum_{k=0}^{n-1}(-1)^{k}\binom{n-1}{k}\int_{0}^{\infty}\ln x\,\, e^{-\left(k+1\right)x}dx,\label{1}
\end{align}
where we have used the substitution $x=e^{y}$ and then applied the Binomial Theorem. Now, for the derivative of the Gamma function at $x=1$, we have
\begin{equation}
\Gamma^{\prime}(1)=\int_{0}^{\infty}e^{-t}\ln t\,\, dt=-\gamma, \label{2}
\end{equation}
where $\gamma=0.57721\cdots$ is the Euler--Mascheroni constant (see, e.g., \cite{B}). Therefore,
\begin{align}
E[Y]=&\left(-n\right)\sum_{k=0}^{n-1}\binom{n-1}{k}\frac{(-1)^{k}}{k+1}\left(\gamma+\ln \left(k+1\right)\right). \label{K}
\end{align}
By invoking the well-known combinatorial identity
\begin{align}
\frac{1}{k+1}\binom{n-1}{k}=\frac{1}{n}\binom{n}{k+1}\label{3}
\end{align}
(see, e.g., \cite{KPG}), relation (\ref{K}) yields
\begin{align*}
E[Y]=&\sum_{j=1}^{n}\left(-1\right)^{j}\binom{n}{j}\left(\gamma+\ln j\right).
\end{align*}
Again, by the Binomial Theorem,
\begin{equation}
\sum_{j=1}^{n}\left(-1\right)^{j}\binom{n}{j}=-1.\label{btm}
\end{equation}
Hence,
\begin{align}
E[Y]=-\gamma+\sum_{j=1}^{n}\left(-1\right)^{j}\binom{n}{j}\ln j. \label{AL}
\end{align}
In view of relation (\ref{cp}), relation (\ref{AL}) yields
\begin{align}
E[Y]=-nc_{n}-\gamma. \label{-1}
\end{align}

Proceeding similarly for the second moment of the random variable $Y$, we obtain
\begin{align}
E[Y^{2}]=&n\int_{-\infty}^{\infty}y^{2}\,e^{y}e^{-e^{y}}\left(1-e^{-e^{y}}\right)^{n-1}dy \nonumber\\
=&n\sum_{k=0}^{n-1}(-1)^{k}\binom{n-1}{k}\int_{0}^{\infty}\left(\ln x\right)^{2}\,\, e^{-\left(k+1\right)x}dx\nonumber\\
=&\sum_{j=1}^{n}\left(-1\right)^{j}\binom{n}{j}\left[\frac{\pi^{2}}{6}+\left(\gamma+\ln j\right)^{2}\right],
\label{4}
\end{align}
where we have also used (\ref{3}), together with the fact that for the second derivative of the Gamma function at $x=1$ we have
\begin{equation}
\Gamma^{\prime\prime}(1)=\int_{0}^{\infty}e^{-x}\left(\ln x\right)^{2}dx=\gamma^{2}+\frac{\pi^{2}}{6}, \label{5}
\end{equation}
(see, e.g., \cite{B}). Next, we expand the square inside the brackets in relation (\ref{4}) and use (\ref{btm}) and (\ref{-1}). Finally, in view of (\ref{cp}) and (\ref{wp}), one arrives at
\begin{align}
E[Y^{2}]=&\gamma^{2}+\frac{\pi^{2}}{6}+2\gamma\, n\,c_{n}+nw_{n}.\label{11}
\end{align}
Having relations (\ref{-1}) and (\ref{11}), the desired result, namely inequality (\ref{v}), follows immediately, since for the strictly positive variance of the random variable $Y$ we have
\begin{equation}
V[Y]=E[Y^{2}]-E[Y]^{2}=\frac{\pi^{2}}{6} - n^{2}c^{2}_{n}+n w_{n},\,\,\,n\in \mathbb{N}.
\label{result}
\end{equation}

\begin{remark}
For the sake of completeness, we briefly mention the case where $n\rightarrow\infty$. In this case, Flajolet and Sedgewick developed a general method for the asymptotic estimation of high-order differences of a fixed numerical sequence $\left\{f_{j}\right\}$ of the form
\begin{equation*}
D_{n}\left[f\right]:=\sum_{j=0}^{n}\binom{n}{j}\left(-1\right)^{j}f_{j}.
\end{equation*}
\end{remark}
We refer the interested reader to \cite{F} and the references therein for details. In \cite{F}, it was proved that as $n\rightarrow \infty$:
\begin{equation}
  \sum_{k=1}^{n}\left(-1\right)^{k}\binom{n}{k}\ln k = \ln\left(\ln n\right)+\gamma+\frac{\gamma}{\ln n}-\frac{\gamma^{2}+\frac{\pi^{2}}{6}}{2\left(\ln n\right)^{2}}+O\left(\frac{1}{\left(\ln n\right)^{3}}\right).\label{fs1}
\end{equation}
By mimicking the method of \cite{F}, and with some patience and computation (since the authors omitted the details), one obtains, as $n\rightarrow \infty$,
\begin{align}
  \sum_{k=1}^{n}\left(-1\right)^{k}\binom{n}{k}\left(\ln k\right)^{2} =&-\left( \ln\left(\ln n\right)\right)^{2}-2\gamma \ln\left(\ln n\right)+\frac{\pi^{2}}{6}-\gamma^{2}-2\gamma \frac{\ln \left(\ln n\right)}{\ln n} \nonumber \\
  &+\frac{\left(\gamma^{2}+\frac{\pi^{2}}{6}\right)\ln \left(\ln n\right)}{\left(\ln n\right)^{2}}
-\frac{2\gamma^{2}}{\ln n}-\frac{\gamma^{2}-\frac{\pi^{2}}{6}}{\left(\ln n\right)^{2}}
  +O\left(\frac{\ln\left(\ln n\right)}{\left(\ln n\right)^{3}}\right).\label{fs2}
\end{align}
By invoking relations (\ref{cp}), (\ref{wp}), and (\ref{result}) in (\ref{fs1}) and (\ref{fs2}), we get
\begin{equation}\label{inf}
 \frac{\pi^{2}}{6} - n^{2}c^{2}_{n}+n w_{n}=o(1),\,\,\,n\rightarrow \infty.
\end{equation}
Hence, for the variance of the random variable $Y$, we have
\begin{equation*}
   \lim_{n\rightarrow \infty}V\left[Y(n)\right]=0.
\end{equation*}
This result arises naturally, since as the number of coupons $N$ and the number of players $n$ both become infinite, one expects that
\begin{equation*}
V\left[M_{N(n)}\right] =o(1).
\end{equation*}
We close with the following remark. In \cite{D}, it was conjectured that the sequence $\frac{\pi^{2}}{6} +n w_{n}- n^{2}c^{2}_{n}$ is decreasing in $n$. This question remains open.\\

\textbf{Acknowledgements.}  This work has been partially supported by project MIS 5154714 of the National Recovery and Resilience Plan Greece 2.0 funded
by the European Union under the NextGenerationEU Program.

\end{document}